# Reversible Digital Filters Total Parametric Sensitivity Optimization using Non-canonical Hypercomplex Number Systems

**Yakiv O. Kalinovsky,** Dr.Sc., Senior Researcher,
Institute for Information Recording National Academy of Science of Ukraine, Kyiv, Shpaka str. 2, 03113, Ukraine, E-mail: kalinovsky@i.ua

**Yuliya E. Boyarinova,** PhD, Associate Professor,
National Technical University of Ukraine "KPI", Kyiv, Peremogy av. 37, 03056, Ukraine,
E-mail: ub@ua.fm

**Iana V. Khitsko,** Junior Researcher,
National Technical University of Ukraine "KPI",Kyiv, Peremogy av. 37, 03056, Ukraine, E-mail: yannuary@yandex.ua

**Abstract**

Digital filter construction method, which is optimal by parametric sensitivity, based on using of non-canonical hypercomplex number systems is proposed and investigated. It is shown that the use of non-canonical hypercomplex number system with greater number of non-zero structure constants in multiplication table can significantly improve the sensitivity of the digital filter.

## 1. Introduction and problem statement

The general approach of using hypercomplex number systems in the construction of the amplitude-frequency characteristics and filter sensitivity calculation is already described in [1-10]. Canonical hypercomplex number systems (HNS) are mainly used in these studies. In this article the attempt was made to develop the methodology for digital filters synthesis in more complex HNS, which contain larger amounts of non-zero structural units in their multiplication tables. As it is shown, this approach allows us to synthesize the reversible digital filters with better parametric sensitivity performance.

## 2. Equivalenting of digital filters with real coefficients and hypercomplex coefficients

In this article the method of synthesis filter structure by converting the digital reversible filter of $n$ order with real coefficients to the digital reversible filter of the first order with hypercomplex coefficients will be used, described in detail in [4].

Consider a digital filter of order 3 with real coefficients, which frequency response is:

$$H_R = \frac{\phi_3 z^{-3} + \phi_2 z^{-2} + \phi_1 z^{-1} + \phi_0}{\varphi_3 z^{-3} + \varphi_2 z^{-2} + \varphi_1 z^{-1} + 1}. \qquad (1)$$

Then, following the procedure described in [4], it can be converted to the digital filter of the first order with a transfer function:

$$H_\Gamma = \frac{A + Bz^{-1}}{1 + Cz^{-1}} = \frac{(A + Bz^{-1}) \cdot \overline{(1 + Cz^{-1})}}{N(1 + Cz^{-1})}, \quad (2)$$

but with hypercomplex coefficients $A, B, C$ which belong to some HNS of dimension 3. In (2) conjugate and norm $N$ are defined in accordance with the formulas defined for HNS, which is used.

Lets consider non-canonical HNS of dimension 3 with multiplication table:

| $\Gamma(e,3)$ | $e_1$ | $e_2$ | $e_3$ |
|---|---|---|---|
| $e_1$ | $e_1$ | $e_2$ | $e_3$ |
| $e_2$ | $e_2$ | $-e_1 + e_3$ | $-2e_2$ |
| $e_3$ | $e_3$ | $-2e_2$ | $2e_1 - e_3$ |

(3)

In the given table there are 4 non-canonical non-zero constants.

HNS $\Gamma(e,3)$ with multiplication table (3) is isomorphic to HNS $R \oplus C$ with the following multiplication table:

| $R \oplus C$ | $E_1$ | $E_2$ | $E_3$ |
|---|---|---|---|
| $E_1$ | $E_1$ | 0 | 0 |
| $E_2$ | 0 | $E_2$ | $E_3$ |
| $E_3$ | 0 | $E_3$ | $-E_2$ |

(4)

As it is shown from the systems (3) and (4) comparation, calculatons in the last one are much easier. This fact can be used to enhance the digital filter performance.

Then the coefficients of the transfer function $H_\Gamma$ are of the form:

$A = a_1 e_1 + a_2 e_2 + a_3 e_3$, $B = b_1 e_1 + b_2 e_2 + b_3 e_3$, $C = c_1 e_1 + c_2 e_2 + c_3 e_3$, $A, B, C \in \Gamma(e,3)$.

When substituted to (2) and transformed we obtain the first-order digital filter transfer function with hypercomplex coefficients in the $\Gamma(e,3)$:

$$H_\Gamma = \frac{a_1 + \dfrac{K}{z} + \dfrac{M}{z^2} + \dfrac{L}{z^3}}{1 + \dfrac{T}{z} + \dfrac{P}{z^2} + \dfrac{Q}{z^3}}, \quad (5)$$

where $K = a_2 c_2 - a_3 c_3 - 3a_1 c_3 + 2a_1 c_1 + b_1$;

$$M = -2b_3c_3 + c_2a_2c_3 + c_2a_2c_1 - 3a_1c_1c_3 + c_2b_2 - 2a_3c_1c_3 + 4a_3c_3^2 + 2a_1c_2^2 + a_1c_1^2 + 2a_3c_2^2 - 3b_1c_3 + 2a_1c_3^2 + 2b_1c_1$$

;

$$L = c_2b_2c_3 + b_1c_1^2 - 2b_3c_1c_3 + c_2b_2c_1 + 2b_1c_2^2 - 3b_1c_1c_3 + 2b_1c_3^2 + 2b_3c_2^2 + 4b_3c_3^2;$$

$$T = 3c_1 - 3c_2; \quad P = -6c_1c_3 + 3c_2^2 + 3c_1^2; \quad Q = 3c_1c_2^2 + 3c_2^2c_3 + c_1^3 - 3c_1^2c_3 + 4c_3^3.$$

Consider a specific example of a third-order filter with real coefficients and the transfer function[1]:

$$H = \frac{0.287589 + 0.6888683 \cdot z^{-1} + 0.6888683 \cdot z^{-2} + 0.287589 \cdot z^{-3}}{1 + 0.418204 \cdot z^{-1} + 0.473048 \cdot z^{-2} + 0.061292 \cdot z^{-3}}.$$

Consider the process of obtaining the coefficient values $a_1, a_2, a_3, b_1, b_2, b_3, c_1, c_2, c_3$.

Using the method described above, we equate the coefficients of the denominator with the same $z^{-i}$, we find the values of hypercomplex coefficients $a_1, a_2, a_3, b_1, b_2, b_3, c_1, c_2, c_3$.

The following system is obtained:

$$\begin{cases} T = 3c_1 - 3c_2 = 0.418204; \\ P = -6c_1c_3 + 3c_2^2 + 3c_1^2 = 0.473048; \\ Q = 3c_1c_2^2 + 3c_2^2c_3 + c_1^3 - 3c_1^2c_3 + 4c_3^3 = 0.061292, \end{cases} \quad (6)$$

whence:

$c_1 = 0.1403252267$,

$c_2 = -0.3718209092$,

$c_3 = 0.0009238933689$.

Substitute these values to the transfer function (5) numerator and equate to the corresponding coefficients of the transfer function of the real filter:

$$\begin{cases} a_1 = 0.287589 \\ 0.2778787733a_1 - 0.3718209092a_2 - 0.001847786738a_3 + b_1 = 0.6888683 \\ -0.001847786738b_3 - 0.05251937625a_2 + 0.2958055168a_1 - 0.3718209092b_2 + \\ + 0.2762457002a_3 + 0.2778787733b_1 = 0.6888683 \\ -0.05251937625b_2 + 0.2958055168b_1 + 0.2762457002b_3 = 0.287589 \end{cases}$$

(7)

Express $a_1, a_2, b_1, b_3$ by $a_3, b_2$:

$$\begin{aligned} a_1 &= 0.287589 \\ a_2 &= 8.446312201 - 5.370129737a_3 + 7.221392887b_2 \\ b_1 &= 3.749468903 - 1.994878735a_3 + 2.685064869b_2 \\ b_3 &= -2.973890946 + 2.136127855a_3 - 2.68506487b_2 \end{aligned} \quad (8)$$

Thus, we obtain the filter parameters as a function of variables $a_3, b_2$. Digital filter parametric sensitivity is the sensitivity of the digital filter transfer function $|H(w)|$ to the coefficients variations of the filter transfer function. Parametric sensitivity function allows us to analyze the impact of coefficients error on the output signal. For the filters with hypercomplex coefficients research we need to consider the possible cumulative error for each of the transfer function coefficients.

Total parametric sensitivity of a first-order filter with hypercomplex coefficients in HNS of dimension 3 is defined by the formula:

$$RCS = \left| \sum_{i=1}^{9} \frac{\alpha_i}{|H|} \cdot \frac{\partial |H|}{\alpha_i} \right|, \qquad (9)$$

where $\alpha_1 = a_1$, $\alpha_2 = a_2$, $\alpha_3 = a_3$, $\alpha_4 = b_1$, $\alpha_5 = b_2$, $\alpha_6 = b_3$, $\alpha_7 = c_1$, $\alpha_8 = c_2$, $\alpha_9 = c_3$.

Parameters $a_3, b_2$ can be of any value in system (8). Suppose, that $a_3 = b_2 = 0$. Then, total parametric sensitivity of the filter with hypercomplex coefficients will be built by the formula (9). Its graph is represented on Fig. 1.

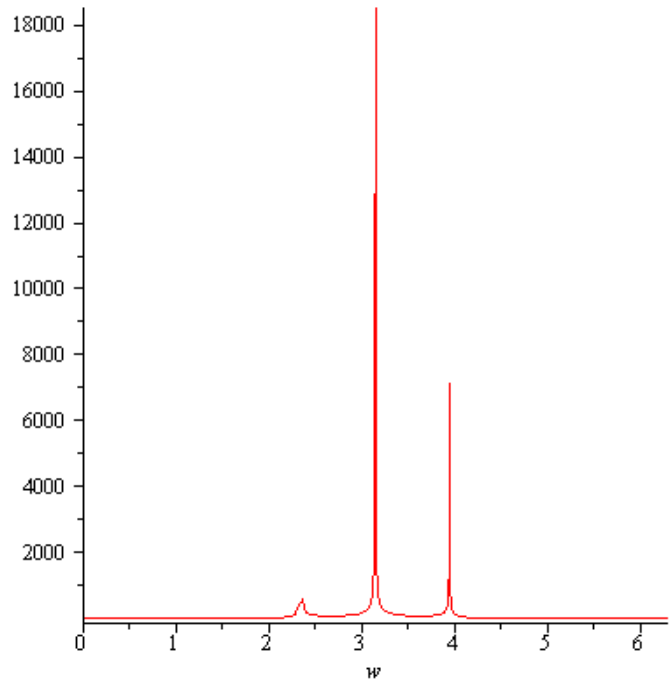

Fig. 1. Total parametric sensitivity of the filter with hypercomplex coefficients.

The total parametric sensitivity of the filter with hypercomplex coefficients to the total parametric sensitivity of the filter with real coefficients ratio plot is presented in Figure 2.

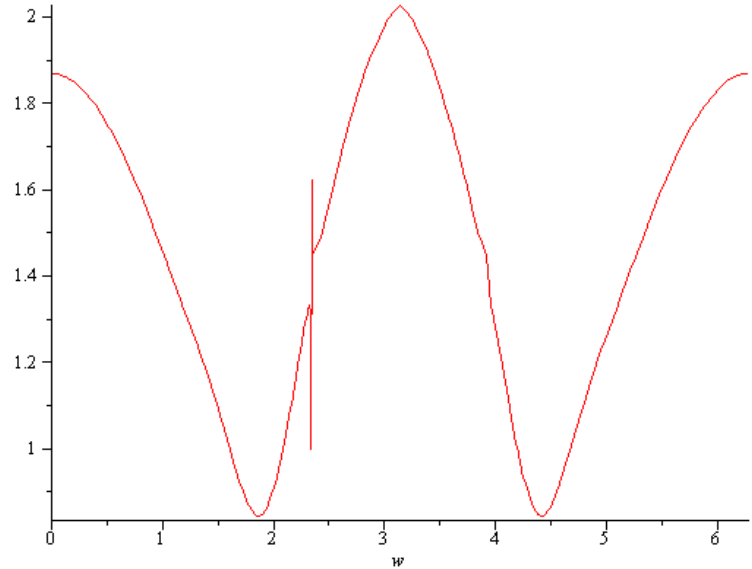

Fig. 2. Ratio of the total parametric sensitivity of the filter with hypercomplex coefficients to the total parametric sensitivity of the filter with real coefficients.

As it can be seen, in this case, the hypercomplex filter sensitivity is much higher than one of a real filter.

## 3. Parametric sensitivity optimization

As it can be seen from (8), the filter parameters $a_3, b_2$ may take different values without altering the transfer function. This fact can be used to optimize the filter parametric sensitivity.

Lets make results optimization. You must choose such values $a_3, b_2$ to satisfy conditions (7) and at the same time to optimize a certain criterion. Calculate total sensitivity function to build this criterion, expressing all its components via $a_3, b_2$.

$RCS = 2.87589 \cdot 10^5 |(z(z^2 + 0.2778787733z + 0.2958055168))/(6.8886829998 \cdot 10^5 z + 6.888683 \cdot 10^5 z^2 +$
$+ 2.87589 \cdot 10^5 z^2 + 2.875889996 \cdot 10^5 - 0.000043zb_2 - 0.0001za_3 - 0.001z^2 a_3 - 0.0004b_2)| +$
$+ 3.718209092 \cdot 10^5 (8.4463312201 - 5.370129737 a_3 + 7.221392887 b_2) | ((z^3 + 0.418204z^2 +$
$+ 0.4730479999z + 0.061292)z)/((z^2 + 0.2769548799z + 0.4339283669)(6.888682998 \cdot 10^5 z +$

$+ 6.888683 \cdot 10^5 z^2 + +2.87589 \cdot 10^5 z^3 + 2.875889996 \cdot 10^5 - 0.000043 zb_2 - 0.0001za_3 - 0.001z^2 a_3 -$
$- 0.0004b_2))| + 2 \cdot 10^6 a_3 |z(0.0009238933\ 689z - 0.1381228502))/(6.888683 \cdot 10^5 z + 6.888683 \cdot 10^5 z^2 +$
$+ 2.87589 \cdot 10^5 z^3 + 2.875889 \cdot 10^5 - 0.000043 zb_2 - 0.0001za_3 - 0.001z^2 a_3 - 0.0004b_2)| + 10^6 (3.749468903$
$- 1.994878735\ a_3 + 2.685064869\ b_2)|z^2 + 0.2778787733\ z + 0.2958055168)/(6.888683 \cdot 10^5 z +$

$+ 0.6888683 \cdot 10^5 z^2 + 2.87589 \cdot 10^5 z^3 + 2.87589 \cdot 10^5 - 0.000043 zb_2 - 0.0001za_3 - 0.001z^2 a_3 - 0.0004b_2)| +$
$+ 3.718209092 \cdot 10^5 b_2 |(z^3 + 0.418204\ z^2 + 0.473047999\ z + 0.061292)/((z^2 + 0.2769548799\ z +$
$+ 0.4339283669)(6.888682998 \cdot 10^5 z^2 + 6.888683 \cdot 10^5 z^2 + 2.87589 \cdot 10^5 z^3 + 2.875889996 \cdot 10^5 -$
$- 0.000043\ zb_2 - 0.0001\ za_3 - 0.001z^2 a_3 - 0.0004b_2))| + 2 \cdot 10^6 (-2.973890946 + 2.136127855\ a_3 -$

$$- 2.68506487b_2 )|(0.0009238933689z + 0.1381228502)/(6.888683 \cdot 10^5 z + 6.888683 \cdot 10^5 z^2 +$$
$$+ 2.87589 \cdot 10^5 z^3 + 2.87589 \cdot 10^5 - 0.000043zb_2 - 0.0001za_3 - 0.001z^2 a_3 - 0.0004b_2 )| + 1.403252267 \cdot 10^5 \cdot$$
$$\cdot |(0.20109286z + 0.8080145342z^2 + 0.3967868515z^3 - 0.07118465934 + 2.685064869z^4 b_2 -$$
$$- 1.994878735z^4 a_3 + 1.502167918z^3 b_2 - 1.39254783z^3 a_3 + 0.3439826239zb_2 + 1.428773891z^2 b_2 -$$

$$- 1.177148979z^2 a_3 - 0.3863640851za_3 - 0.287589z^5 + 2.371732304 - 0.0342181936a_3 + 0.02324579b_2 )$$
$$/((6.888683 \cdot 10^5 z + 6.888683 \cdot 10^5 z^2 + 2.87589 \cdot 10^5 z^3 + 2.87589 \cdot 10^5 - 0.000043zb_2 - 0.0001za_3 -$$
$$- 0.001z^2 a_3 - 0.0004b_2 )(z^3 + 0.488204z^2 + 0.473048z + 0.061292))| - 3.71820992 \cdot 10^5 |((3.899734973z +$$
$$+ 2.553110588z^2 + 0.9985997616 + 8.446312201z^3 - 5.370129737z^3 a_3 + 3z^2 b_2 - 0.6453745599a_3 +$$

$$+ 0.439283673 b_2 )(z^3 + 8.446312201 z^3 - 5.370129737 z^3 a_3 + 3z^2 b_2 - 0.6453745599 a_3 + 0.439283673 b_2 ) \cdot$$
$$\cdot (z^3 + 0.418204 z^2 + 0.473048 z + 0.061292 ))/((z^2 + 0.2769548799 z + 0.4339283669 )^2 (6.888682998 \cdot 10^5 z$$
$$+ 6.888683 \cdot 10^5 z^2 + 2.87589 \cdot 10^5 z^3 + 2.875889996 \cdot 10^5 - 0.000043 zb_2 - 0.0001za_3 - 0.001z^2 a_3 -$$
$$- 0.0004 b_2 ))| + 1847.786738 |(- 0.4416591322 z - 1.609719097 z^2 - 1.098032331 z^3 - 0.0741899487 -$$

$$- 2.685064868 z^4 b_2 + 1.299719855 z^4 a_3 - 1.502167918 z^3 b_2 + 0.3697478256 z^3 a_3 - z^5 a_3 - 0.3439826232 zb_2$$
$$- 1.42877389 z^2 b_2 + 0.8033734483 z^2 a_3 + 0.1641200202 za_3 - 3.306608615 z^4 + 0.0076218562 a_3 -$$
$$- 0.0232457908 b_2 )/((6.8886823 \cdot 10^5 z + 6.888683 \cdot 10^5 z^2 + 2.87589 \cdot 10^5 z^3 + 2.875889 \cdot 10^5 - 0.000043 zb_2$$
$$- 0.0001za_3 - 0.001z^2 a_3 - 0.0004b_2 )(z^3 + 0.418204 z^2 + 0.473048 z + 0.06129 ))|$$

Since the sensitivity function is positive on the whole interval $\omega = 0..2\pi$, the optimality criterion can take the parametric sensitivity sum for a certain set of values $\omega$, with parameters values $a_3, b_2$. We select the 33 evenly distributed points on the interval $\{0..2\pi\}$ and calculate the function values at each point, given the fact that $z = \sin(\omega) + i \cdot \cos(\omega)$. Then construct the optimality criterion $S_{RCS}(\omega, a_3, b_2)$ which must be minimized.

It is impossible to represent the function $S_{RCS}(\omega, a_3, b_2)$ in this article, since it is too cumbersome. It was also unsuccessful to apply a gradient optimization method, such as function differentiation by the components $a_3, b_2$ is very complex. Therefore, its optimization is an independent honest task.

To prove the described method efficiency of digital filter synthesis is sufficient to find the approximate optimum, which is possible to be done by the construction of a function three-dimensional graph, which used procedures of analytical calculations MAPLE. At the same time it can be a multistage procedure: first select the wide scope of the search, then it narrows. Accordingly, on Fig. 3. a wide-range search area is presented, on Fig. 4. – narrowed one.

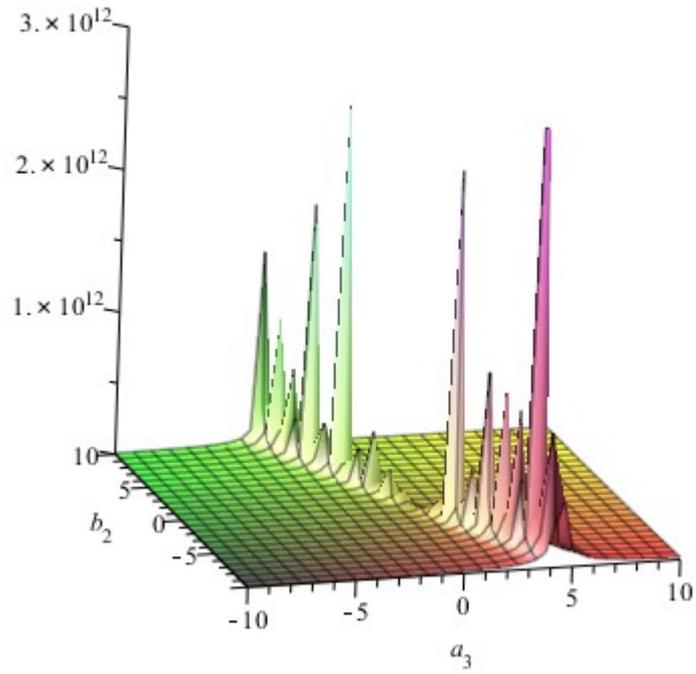

Fig. 3. $S_{RCS}(\omega, a_3, b_2)$ plot wide range search area; $a_3 \in \{-10..10\}, b_2 \in \{-10..10\}$.

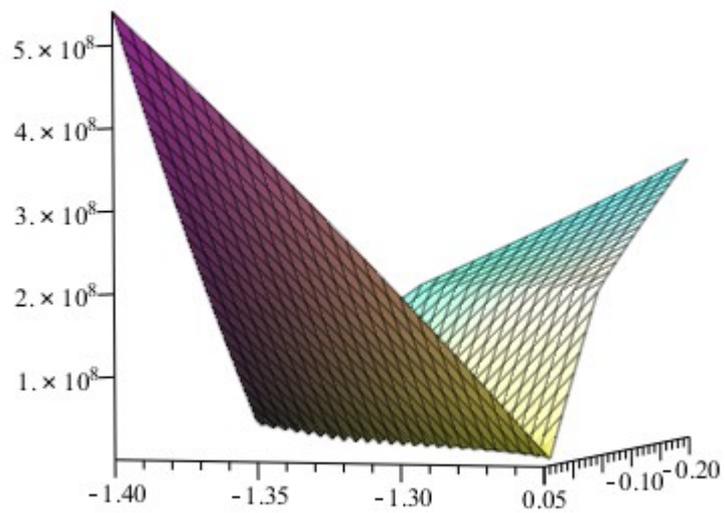

Fig. 4. $S_{RCS}(\omega, a_3, b_2)$ plot for narrow range search area; $a_3 \in \{-0.25..0.05\}, b_2 \in \{-1.4..-1.25\}$

Figure 5 is showing the parametric sensitivity changes graph near one of the obtained local minimum with $a_3 = -0.2316615, b_2 = -1.2783899677$.

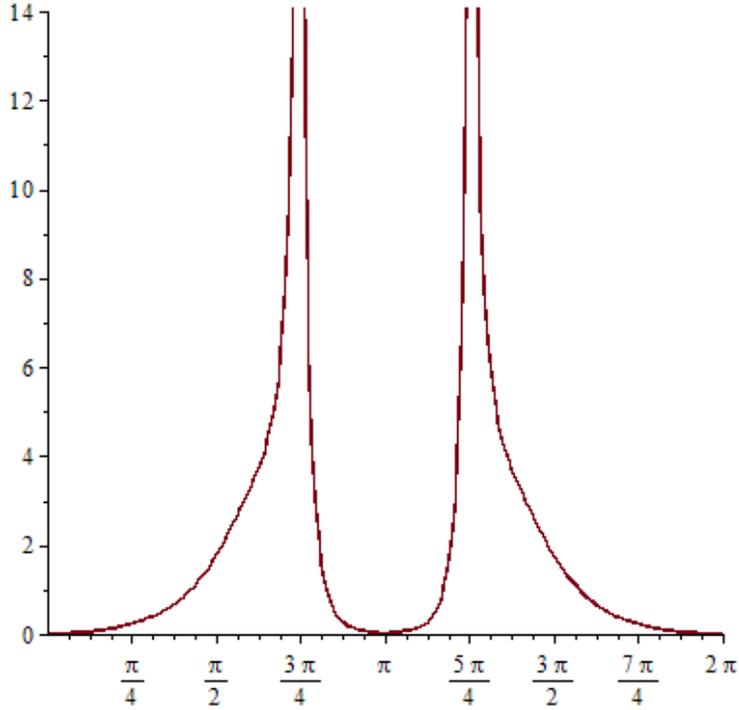

Fig. 5. Total parametric sensitivity of the filter with hypercomplex coefficients after optimization.

The sensitivity of the filter with hypercomplex coefficients in $\Gamma(e,3)$ system to the sensitivity of the filter with real coefficients ratio is shown on Fig. 6, which shows that the hypercomplex filter sensitivity is lower than the one of the real filter.

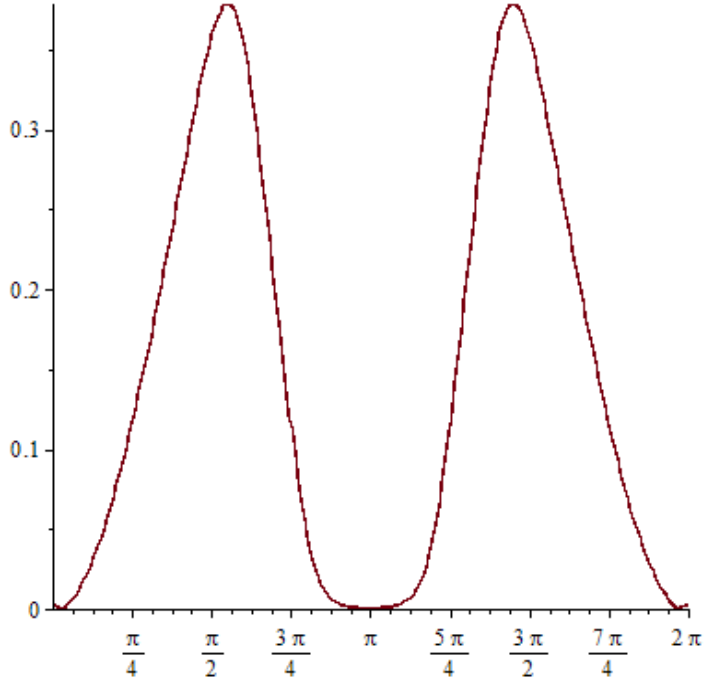

Fig. 6. The sensitivity of the filter with hypercomplex coefficients in $\Gamma(e,3)$ to the sensitivity of the filter with real coefficients ratio.

**4. Conclusions**

Thus, we have shown that using the non-canonical HNS allows to reduce the total parametric sensitivity of the digital filter. At the same time, precise optimization of the target function $S_{RCS}(\omega, a_3, b_2)$ requires additional research.